\documentclass[leqno,11pt]{amsart}

\usepackage[T1]{fontenc}
\usepackage[utf8]{inputenc}
\usepackage[a4paper,vmargin={3cm,3cm},hmargin={3cm,3cm}]{geometry}
\usepackage[english]{babel}
\usepackage{amsmath,amsfonts,amssymb,mathtools}
\usepackage{hyperref}
\usepackage[nopatch=showhyphens]{microtype}
\numberwithin{equation}{section}

\newtheorem{thm}{Theorem}[section]
\newtheorem{lem}[thm]{Lemma}

\theoremstyle{remark}
\newtheorem{rem}[thm]{Remark}

\title[Improvement of effective Erd\H{o}s--Wintner theorem]
{Improvement of effective Erd\H{o}s--Wintner theorem\\for Zeckendorf expansions}

\author{Johann Verwee}

\begin{document}

\begin{abstract}
In their effective Erd\H{o}s--Wintner theorem for Zeckendorf expansions, Drmota and the author obtained a uniform Kolmogorov bound with an error term that involves
\( T\sum_{j>L-2h} |f(F_j)| \),
hence assumes absolute convergence of the linear tail \(\sum_j f(F_j)\).
We show that the absolute convergence can be removed. The key is to group the transfer matrices in pairs and work at second order on the logarithm of the product, after extracting the common linear phase along the dominant direction. This yields a \emph{quadratic} tail
\( T^2\sum_{j>L-2h} f(F_j)^2 \),
or, in a flexible variant, the split tail
\( T\sum_{|f(F_j)|>1/T} |f(F_j)| + T^2 \sum_{|f(F_j)|\leqslant 1/T} f(F_j)^2 \).
Either form requires only \(\sum f(F_j)^2<\infty\).
\end{abstract}


\maketitle

\section{Main result}
\begin{thm}\label{thm:improve}
Let \(F_0=0\), \(F_1=1\), \(F_{j+1}=F_j+F_{j-1}\), and write \(G_k:=F_{k+2}\).
Let \(f\) be a real Zeckendorf-additive function, i.e.\
\(f(n)=\sum_{j\ge2} f(\epsilon_j F_j)\) for the Zeckendorf digits \(\epsilon_j\in\{0,1\}\) (no two consecutive \(1\)'s).
For \(k\ge0\) set \(H_k(t):=\sum_{0\leqslant n<G_k}e^{it f(n)}\) and \(\Phi_k(t):=H_k(t)/H_k(0)\),
and let \(\Phi(t)=\lim_{k\to\infty}\Phi_k(t)\) be the characteristic function of the limit law \(F\).
Let \(Q_F\) denote the concentration function of \(F\).
Assume \(\sum_{j\ge0} f(F_j)^2<\infty\).
Then, for \(N\) large and parameters \(h,L\) chosen as in \cite[Th.\,8]{DV},
\begin{equation}\label{eq:main}
\|F_N-F\|_\infty
\ \ll\ Q_F(1/T)\ +\ \frac{\log N}{T}\ +\ T^2\sum_{j>L-2h} f(F_j)^2,
\end{equation}
uniformly for \(T\in(0,1]\). In particular, absolute convergence of \(\sum_j f(F_j)\) is not required.

Moreover, one has the flexible split form
\begin{equation}\label{eq:split}
\|F_N-F\|_\infty
\ \ll\ Q_F(1/T)\ +\ \frac{\log N}{T}\ +\ T\sum_{\;|f(F_j)|>1/T}|f(F_j)|
\ +\ T^2\sum_{\;|f(F_j)|\leqslant 1/T} f(F_j)^2.
\end{equation}
For fixed \(T\), the set \(\{j:\ |f(F_j)|>1/T\}\) is finite under \(\ell^2\), so \eqref{eq:split} also needs only \(\ell^2\).
\end{thm}

\section{Preliminaries and notation}
Let \(F_0=0\), \(F_1=1\), \(F_{j+1}=F_j+F_{j-1}\). Every \(n\ge0\) has a unique Zeckendorf expansion
\(n=\sum_{j\ge2}\epsilon_j F_j\) with \(\epsilon_j\in\{0,1\}\) and no two consecutive \(1\)'s \cite{Zeckendorf}.
Set
\[
G_k:=F_{k+2}\qquad(k\ge0).
\]
Fix values \(f(F_j)\) for \(j\ge2\) and set \(f(0):=0\). Extend \(f\) additively by
\[
f(n)=\sum_{j\ge2} \epsilon_j\,f(F_j),
\]
where \(\epsilon_j\in\{0,1\}\) are the Zeckendorf digits.

For \(t\in\mathbb{R}\), define the characteristic sum at height \(k\) by
\begin{equation}\label{eq:def-Hk}
H_k(t):=\sum_{0\leqslant n<G_k} e^{i t f(n)}\qquad(G_k=F_{k+2}),
\end{equation}
and the normalized characteristic function \(\Phi_k(t):=H_k(t)/H_k(0)\).
The limit (when it exists) \(\Phi(t)=\lim_{k\to\infty}\Phi_k(t)\) is the characteristic function of the limiting distribution \(F\).

\subsection*{Scalar recursion and transfer matrices}
A direct conditioning on the last Zeckendorf digit gives, for \(k\geqslant 1\),
\begin{equation}\label{eq:scalar-rec}
H_{k+1}(t)=H_k(t)+e^{i t f(F_{k+2})}\,H_{k-1}(t),
\qquad H_0(t)=1,\ H_1(t)=1+e^{i t f(F_2)}.
\end{equation}
Equivalently, in \(2\times2\) form,
\begin{equation}\label{eq:matrix-rec}
\binom{H_{k+1}(t)}{H_k(t)}
=
A_k(t)\binom{H_k(t)}{H_{k-1}(t)},
\qquad
A_k(t):=\begin{pmatrix}1&1\\ e^{i t f(F_{k+2})}&0\end{pmatrix}.
\end{equation}
At \(t=0\) one gets the constant transfer matrix
\[
A:=A_k(0)=\begin{pmatrix}1&1\\ 1&0\end{pmatrix},
\]
which is diagonalizable over \(\mathbb{R}\) with eigenvalues
\[
\alpha=\tfrac{1+\sqrt5}{2},\qquad
\overline{\alpha}=\tfrac{1-\sqrt5}{2}=-\alpha^{-1},\qquad \alpha^2=\alpha+1.
\]
We write
\begin{equation}\label{eq:def-Delta}
\Delta_k(t):=A_k(t)-A
=\begin{pmatrix}0&0\\ e^{i t f(F_{k+2})}-1&0\end{pmatrix},
\quad\text{so that } A_k(t)=A+\Delta_k(t).
\end{equation}

\subsection*{Basic bounds and adapted norm}
For \(|t|\leqslant 1\), the inequality
\[
|e^{ix}-1|\leqslant \min\{|x|,x^2\}\leqslant |x|+x^2
\]
gives the uniform bound (use the Frobenius norm, which dominates any fixed operator norm)
\begin{equation}\label{eq:Delta-bound}
\|\Delta_k(t)\|\ \le\ \|\Delta_k(t)\|_{\mathrm{F}}
= |e^{i t f(F_{k+2})}-1|
\ \ll\ |t|\,|f(F_{k+2})|+t^2 f(F_{k+2})^2.
\end{equation}
Let \(P\) diagonalize \(A\): \(P^{-1}AP=\mathrm{diag}(\alpha,\overline{\alpha})\).
Define the \(A\)-adapted norm
\begin{equation}\label{eq:adapted-norm}
\|M\|_{A}:=\|P^{-1}MP\|_2,
\end{equation}
so that \(\|A\|_{A}=\alpha\), \(\|A^2\|_{A}=\alpha^2\), and
\begin{equation}\label{eq:Delta-adapted}
\|P^{-1}\Delta_k(t)P\|_{A}\ \ll_{A}\ \|\Delta_k(t)\|
\ \ll\ |t|\,|f(F_{k+2})|+t^2 f(F_{k+2})^2.
\end{equation}

\subsection*{Block factorization and cancellation of linear terms}
Group the one–step matrices two by two:
\begin{equation}\label{eq:def-Bk}
\mathcal{B}_k(t):=A_{2k+1}(t)\,A_{2k}(t)
= A^2 + A\,\Delta_{2k}(t) + \Delta_{2k+1}(t)\,A + \Delta_{2k+1}(t)\Delta_{2k}(t).
\end{equation}
With \eqref{eq:def-Delta}, one has the exact identity
\[
\Delta_{2k+1}(t)\Delta_{2k}(t)= 0
\quad\text{since}\quad
\begin{pmatrix}0&0\\ *&0\end{pmatrix}\begin{pmatrix}0&0\\ *&0\end{pmatrix}=0.
\]
Conjugating by \(P\) and writing tildes for conjugated matrices gives
\begin{equation}\label{eq:tilde-Bk}
\widetilde{\mathcal{B}}_k(t)
:=P^{-1}\mathcal{B}_k(t)P
=
D^2 \;+\; D\,\widetilde{\Delta}_{2k}(t)\;+\;\widetilde{\Delta}_{2k+1}(t)\,D.
\end{equation}

\begin{lem}[Explicit linear term in the dominant coordinate]\label{lem:lin-phase}
Let \(E_{21}\) denote the \(2\times2\) matrix with a single \(1\) in position \((2,1)\) and zeros elsewhere.
For $m\in\mathbb{Z}$ set $\delta_m(t):=e^{it\,f(F_{m+2})}-1$, $D:=\mathrm{diag}(\alpha,\overline{\alpha})$,
and $S:=P^{-1}E_{21}P$. Then
\begin{equation}\label{eq:tilde-Bk-exact}
\widetilde{\mathcal{B}}_k(t)=D^2+\delta_{2k}(t)\,DS+\delta_{2k+1}(t)\,SD.
\end{equation}
In particular, for the $(1,1)$ entry,
\begin{equation}\label{eq:lin-11-exact}
\big(\widetilde{\mathcal{B}}_k(t)\big)_{11}
=\alpha^2+\frac{\alpha}{\sqrt5}\,\big(\delta_{2k}(t)+\delta_{2k+1}(t)\big).
\end{equation}
Moreover, writing $\delta_m(t)=it\,f(F_{m+2})+r_m(t)$ with $|r_m(t)|\ll t^2 f(F_{m+2})^2$ for $|t|\leqslant 1$,
one gets
\begin{equation}\label{eq:lin-11-expanded}
\big(\widetilde{\mathcal{B}}_k(t)\big)_{11}
=\alpha^2+\frac{i\alpha}{\sqrt5}\,t\,(f(F_{2k+2})+f(F_{2k+3}))
\;+\;O_A\big(t^2(f(F_{2k+2})^2+f(F_{2k+3})^2)\big).
\end{equation}
\end{lem}

\bigskip

\begin{proof}
By definition $\Delta_m(t)=A_m(t)-A=\delta_m(t)\,E_{21}$, hence
\[
\mathcal{B}_k(t)=A_{2k+1}(t)A_{2k}(t)=A^2+A\,\Delta_{2k}(t)+\Delta_{2k+1}(t)\,A,
\]
and, after conjugation by $P$,
\[
\widetilde{\mathcal{B}}_k(t)
= D^2 \;+\; \delta_{2k}(t)\,P^{-1}A E_{21}P\;+\;\delta_{2k+1}(t)\,P^{-1}E_{21}A P.
\]
Let $S:=P^{-1}E_{21}P$. Compute $P^{-1}$ from
\[
P=\begin{pmatrix}\alpha&\overline{\alpha}\\ 1&1\end{pmatrix},
\qquad
\det P=\alpha-\overline{\alpha}=\sqrt5,
\qquad
P^{-1}=\frac1{\sqrt5}\begin{pmatrix}1&-\overline{\alpha}\\ -1&\alpha\end{pmatrix}.
\]
Since $E_{21}P=\bigl(\begin{smallmatrix}0&0\\ \alpha&\overline{\alpha}\end{smallmatrix}\bigr)$, we get
\[
S=P^{-1}E_{21}P
=\frac1{\sqrt5}\begin{pmatrix}1&-\overline{\alpha}\\ -1&\alpha\end{pmatrix}
\begin{pmatrix}0&0\\ \alpha&\overline{\alpha}\end{pmatrix}
=\frac1{\sqrt5}\begin{pmatrix}1&-\overline{\alpha}^{\,2}\\ \alpha^{2}&-1\end{pmatrix},
\]
using $\alpha\overline{\alpha}=-1$ so that $\overline{\alpha}^{\,2}=1-\alpha$ and $\alpha^{2}=\alpha+1$.
As $P^{-1}AP=D$, we have
\[
P^{-1}A E_{21}P=DS
\qquad\text{and}\qquad
P^{-1}E_{21}A P=SD.
\]
Read off their $(1,1)$ entries:
\[
(DS)_{11}=\frac{\alpha}{\sqrt5},\qquad (SD)_{11}=\frac{\alpha}{\sqrt5}.
\]
Plugging this into $\widetilde{\mathcal{B}}_k(t)$ yields
\[
\big(\widetilde{\mathcal{B}}_k(t)\big)_{11}
=\alpha^2+\frac{\alpha}{\sqrt5}\,\delta_{2k}(t)+\frac{\alpha}{\sqrt5}\,\delta_{2k+1}(t),
\]
which is exactly \eqref{eq:lin-11-exact}. Finally, $\delta_m(t)=it\,f(F_{m+2})+O(t^2 f(F_{m+2})^2)$
gives the linear term $\frac{i\alpha}{\sqrt5}\,t\,(f(F_{2k+2})+f(F_{2k+3}))$.
\end{proof}

Let \(E_{11}=\bigl(\begin{smallmatrix}1&0\\[.2mm]0&0\end{smallmatrix}\bigr)\) and set
\[
M_k(t):=D^{-2}\,\widetilde{\mathcal{B}}_k(t)
=I+\delta_{2k}(t)\,U+\delta_{2k+1}(t)\,V,
\qquad U:=D^{-2}DS,\quad V:=D^{-2}SD.
\]
Then \(U_{11}=V_{11}=1/(\alpha\sqrt5)\).
Define the phase parameter
\[
\vartheta_k(t):=\frac{\delta_{2k}(t)+\delta_{2k+1}(t)}{\alpha\sqrt5}
=\frac{i t}{\alpha\sqrt5}\big(f(F_{2k+2})+f(F_{2k+3})\big)+O\big(t^2(f(F_{2k+2})^2+f(F_{2k+3})^2)\big).
\]
Since \(E_{11}^2=E_{11}\), we have \(e^{\pm\vartheta_k E_{11}}=I+(e^{\pm\vartheta_k}-1)E_{11}\).
Define
\[
\begin{aligned}
R_k(t)
&:= e^{-\vartheta_k(t)E_{11}}\,M_k(t) - I \\
&= (e^{-\vartheta_k}-1)E_{11} + \delta_{2k}U + \delta_{2k+1}V \\
&\quad + (e^{-\vartheta_k}-1)E_{11}\bigl(\delta_{2k}U+\delta_{2k+1}V\bigr).
\end{aligned}
\]
\paragraph{Exact (1,1)-entry and quadratic bound.}
Recall that $E_{11}^2=E_{11}$, hence
\[
e^{-\vartheta_k(t) E_{11}} \;=\; I + \bigl(e^{-\vartheta_k(t)}-1\bigr)E_{11}.
\]
With $M_k(t)=I+\delta_{2k}(t)\,U+\delta_{2k+1}(t)\,V$, $U_{11}=V_{11}=\tfrac{1}{\alpha\sqrt5}$, and
$\vartheta_k(t)=\tfrac{\delta_{2k}(t)+\delta_{2k+1}(t)}{\alpha\sqrt5}$, we obtain
\begin{align*}
(R_k(t))_{11}
&=\bigl(e^{-\vartheta_k(t) E_{11}}M_k(t)-I\bigr)_{11}\\
&=\bigl(e^{-\vartheta_k(t)}-1\bigr)\,\bigl(M_k(t)\bigr)_{11}+\bigl(M_k(t)\bigr)_{11}-1\\
&=\bigl(e^{-\vartheta_k(t)}-1\bigr)\Bigl(1+\delta_{2k}(t)\,U_{11}+\delta_{2k+1}(t)\,V_{11}\Bigr)
 \;+\;\delta_{2k}(t)\,U_{11}+\delta_{2k+1}(t)\,V_{11}\\
&=\bigl(e^{-\vartheta_k(t)}-1\bigr)\bigl(1+\vartheta_k(t)\bigr)\;+\;\vartheta_k(t)
\;=\; e^{-\vartheta_k(t)}\bigl(1+\vartheta_k(t)\bigr)-1.
\end{align*}
Since, as shown above,
\[
(R_k(t))_{11}=e^{-\vartheta_k(t)}\bigl(1+\vartheta_k(t)\bigr)-1,
\]
we now expand the scalar exponential at $\vartheta_k(t)=0$:
\begin{align*}
e^{-\vartheta_k(t)}\bigl(1+\vartheta_k(t)\bigr)-1
&=\Bigl(1-\vartheta_k(t)+\tfrac{\vartheta_k(t)^2}{2}+O(|\vartheta_k(t)|^3)\Bigr)\bigl(1+\vartheta_k(t)\bigr)-1\\
&= -\tfrac{1}{2}\,\vartheta_k(t)^2\;+\;O\big(|\vartheta_k(t)|^3\big).
\end{align*}
In particular, the linear term cancels. Using \(|\delta_m(t)|\leqslant |t|\,|f(F_{m+2})|+t^2 f(F_{m+2})^2\) for \(|t|\leqslant 1\),
we have \(|\vartheta_k(t)|\ll_A |\delta_{2k}(t)|+|\delta_{2k+1}(t)|\), hence
\[
(R_k(t))_{11}
=O_A\big(t^2\,(f(F_{2k+2})^2+f(F_{2k+3})^2)\big).
\]

\medskip
\noindent\textbf{A-adapted norm bound (independent of the previous cancellation).}
From the definition
\[
R_k(t)=\bigl(e^{-\vartheta_k(t)}-1\bigr)E_{11}\;+\;\delta_{2k}(t)\,U\;+\;\delta_{2k+1}(t)\,V
\;+\;\bigl(e^{-\vartheta_k(t)}-1\bigr)E_{11}\bigl(\delta_{2k}(t)\,U+\delta_{2k+1}(t)\,V\bigr),
\]
and, by the triangle inequality and submultiplicativity of \(\|\cdot\|_A\),
\[
\begin{aligned}
\|R_k(t)\|_{A}
&\leqslant |e^{-\vartheta_k(t)}-1|\,\|E_{11}\|_{A}
   + |\delta_{2k}(t)|\,\|U\|_{A}
   + |\delta_{2k+1}(t)|\,\|V\|_{A} \\
&\quad + |e^{-\vartheta_k(t)}-1|\!\left(
      |\delta_{2k}(t)|\,\|E_{11}U\|_{A}
    + |\delta_{2k+1}(t)|\,\|E_{11}V\|_{A}
    \right).
\end{aligned}
\]
Since \(U=D^{-1}S\) and \(V=D^{-2}SD\), the matrices \(E_{11},U,V\) (and \(E_{11}U,E_{11}V\)) are fixed and depend
only on \(A\). Writing
\(\|E_{11}\|_{A}\leqslant c_1,\ \|U\|_{A}\leqslant c_2,\ \|V\|_{A}\leqslant c_3,\ \|E_{11}U\|_{A}\leqslant c_4,\ \|E_{11}V\|_{A}\leqslant c_5\),
we get
\[
\|R_k(t)\|_{A}
\leqslant C_A\Big(
 |e^{-\vartheta_k(t)}-1| + |\delta_{2k}(t)| + |\delta_{2k+1}(t)|
 + |e^{-\vartheta_k(t)}-1|\big(|\delta_{2k}(t)|+|\delta_{2k+1}(t)|\big)
\Big),
\]
with \(C_A:=c_1+c_2+c_3+c_4+c_5\). For \(|\vartheta_k(t)|\leqslant 1\) (true for large \(k\) under \(\ell^2\)), \(|e^{-\vartheta_k(t)}-1|\leqslant 2|\vartheta_k(t)|\),
and \(|\vartheta_k(t)|\ll_A |\delta_{2k}(t)|+|\delta_{2k+1}(t)|\). Using
\(|\delta_m(t)|\leqslant |t|\,|f(F_{m+2})|+t^2 f(F_{m+2})^2\) finally yields
\[
\|R_k(t)\|_{A}\ \ll_A\ |t|\,(|f(F_{2k+2})|+|f(F_{2k+3})|)\;+\;t^2\,(f(F_{2k+2})^2+f(F_{2k+3})^2).
\]

\bigskip

Moreover, we have the exact factorization
\begin{equation}\label{eq:paired-factor}
\widetilde{\mathcal{B}}_k(t)
= D^2\,M_k(t) = D^2\,\exp\Big(\vartheta_k(t)\,E_{11}\Big)\,\big(I+R_k(t)\big)
\end{equation}
and the product over \(k\) yields
\begin{equation}\label{eq:paired-product}
\prod_{k=K}^{L}\widetilde{\mathcal{B}}_k(t)
= D^{2(L-K+1)}\exp\Big(\sum_{k=K}^{L}\vartheta_k(t)\,E_{11}\Big)
\cdot \prod_{k=K}^{L}\big(I+R_k(t)\big).
\end{equation}
Here and below we set \(V_k:=f(F_{2k+2})^2+f(F_{2k+3})^2\). Since
\begin{equation}\label{eq:phase-size}
\sum_{k=K}^{L}\vartheta_k(t)
= i\frac{t}{\alpha\sqrt5}\sum_{k=K}^{L}\big(f(F_{2k+2})+f(F_{2k+3})\big)
+ O\Big(t^2\sum_{k=K}^{L} V_k\Big),
\end{equation}
the exponential factor is unimodular up to
\(\exp\big(O(t^2\sum_{k=K}^{L} V_k)\big)\), which is absorbed by the same
quadratic tail bound. Hence this phase cancels in \(|\Phi_N(t)-\Phi(t)|\).

\section{Proof of Theorem \ref{thm:improve}}
Let
\[
\widetilde{\mathcal{P}}_{K,L}(t):=\prod_{k=K}^{L}\big(I+R_k(t)\big).
\]
By the bounds obtained above,
\[
(R_k(t))_{11}=O_A\big(t^2 V_k\big)
\quad\text{and}\quad
\|R_k(t)\|_{A}\ \ll_A\ |t|\big(|f(F_{2k+2})|+|f(F_{2k+3})|\big)+t^2 V_k.
\]
Since \(f(F_j)\to0\) (because \(\sum_j f(F_j)^2<\infty\)), for any fixed \(T\in(0,1]\) one can take \(K\) so large that
\[
\sup_{k\geqslant K}\|R_k(t)\|_{A}\leqslant \tfrac12 \qquad\text{for all } |t|\leqslant T.
\]
Hence the matrix logarithm expansion is valid and yields
\[
\log \widetilde{\mathcal{P}}_{K,L}(t)
=\sum_{k=K}^{L}R_k(t)\ +\ O\Big(\sum_{k=K}^{L}\|R_k(t)\|_{A}^2\Big).
\]
Using \(\big(|f(F_{2k+2})|+|f(F_{2k+3})|\big)^2\leqslant 2V_k\) and \(T\leqslant 1\), we get
\[
\|R_k(t)\|_{A}^2\ \ll_A\ t^2 V_k,
\]
so that, entrywise,
\begin{equation}\label{eq:block-tail-11}
\left|(\widetilde{\mathcal{P}}_{K,L}(t))_{11}-1\right|
\ \ll_{A}\ t^2\sum_{k=K}^{L} V_k
\ =\ t^2\sum_{j=2K+2}^{2L+3} f(F_j)^2.
\end{equation}

Fix \(T\in(0,1]\). By the paired-block factorization and the phase extraction proved above,
for the truncation level \(K\asymp L-2h\) one has
\[
\prod_{k=K}^{L}\widetilde{\mathcal{B}}_k(t)
= D^{2(L-K+1)}\exp\Big(\sum_{k=K}^{L}\vartheta_k(t)\,E_{11}\Big)\cdot \widetilde{\mathcal{P}}_{K,L}(t),
\qquad
\widetilde{\mathcal{P}}_{K,L}(t):=\prod_{k=K}^{L}(I+R_k(t)).
\]
By \eqref{eq:block-tail-11} (together with \((R_k)_{11}=O_A(t^2V_k)\) and
\(\| R_k\|_A^2\ll_A t^2V_k\)), we get the block-tail bound
\[
\left|(\widetilde{\mathcal{P}}_{K,L}(t))_{11}-1\right|
\ \ll_A\ t^2\sum_{k=K}^{L}V_k
\ =\ t^2\sum_{j=2K+2}^{2L+3} f(F_j)^2.
\]
As in \cite[Th.\,8]{DV}, writing \(H_N(t)\) as a product of such paired blocks acting on an initial vector that depends only on \(O(1)\) layers, and cutting at level \(L-2h\), the comparison between the finite product (defining \(\Phi_N\)) and the infinite one (defining \(\Phi\)) yields, uniformly for \(|t|\leqslant T\),
\begin{equation}\label{eq:Phi-tail}
|\Phi_N(t)-\Phi(t)|
 \ \ll_A\ T^2\sum_{j>L-2h} f(F_j)^2.
\end{equation}
(Here we used that the phase factor is unimodular up to \(e^{O(t^2\sum V_k)}\), which is absorbed by the same quadratic tail bound, and then took the supremum over \(|t|\leqslant T\).)

\subsection*{Smoothing, choice of \texorpdfstring{$T$}{TEXT}, and atomic vs.\ non-atomic limits}
We shall use the standard Feller--Esseen smoothing inequality (see \cite[Ch.~XVI]{FellerII} or \cite{EsseenCF})
\begin{equation}\label{eq:Esseen-correct}
\|F_N-F\|_\infty \ \ll\ Q_F(1/T)\ +\ \frac{1}{T}\int_0^T \frac{|\Phi_N(t)-\Phi(t)|}{t}\,dt,
\qquad T\ge 1.
\end{equation}
We split the integral at $t_0:=N^{-1}$ and at $1$. On $[0,t_0]$ we argue as in \cite[Th.\,8]{DV} and obtain
\begin{equation}\label{eq:lf}
\frac{1}{T}\int_0^{t_0} \frac{|\Phi_N(t)-\Phi(t)|}{t}\,dt\ \ll\ \frac{\log N}{T}.
\end{equation}
On $[t_0,1]$ we use the paired-block bound \eqref{eq:Phi-tail} with $T=1$, which yields
\begin{equation}\label{eq:mid}
\frac{1}{T}\int_{t_0}^{1} \frac{|\Phi_N(t)-\Phi(t)|}{t}\,dt
\ \ll\ \frac{1}{T}\int_{t_0}^{1}\frac{dt}{t}\ \sum_{j>L-2h} f(F_j)^2
\ \ll\ \frac{\log N}{T}\sum_{j>L-2h} f(F_j)^2.
\end{equation}
Finally, on $[1,T]$ we use the trivial bound $|\Phi_N|,|\Phi|\le 1$ to get
\begin{equation}\label{eq:hf}
\frac{1}{T}\int_{1}^{T} \frac{|\Phi_N(t)-\Phi(t)|}{t}\,dt\ \ll\ \frac{\log T}{T}.
\end{equation}
Combining \eqref{eq:Esseen-correct}–\eqref{eq:hf} and the tail control \eqref{eq:Phi-tail} gives, for all $T\ge 1$,
\begin{equation}\label{eq:master}
\|F_N-F\|_\infty \ \ll\ Q_F(1/T)\ +\ \frac{\log N}{T}
\ +\ \frac{\log N}{T}\sum_{j>L-2h} f(F_j)^2\ +\ \frac{\log T}{T}.
\end{equation}
Since $\sum_j f(F_j)^2<\infty$ and $L=L(N)\to\infty$, the tail $\sum_{j>L-2h} f(F_j)^2\to0$ as $N\to\infty$.

\medskip
\noindent\textbf{Non-atomic limit.}
By \cite[Prop.~11]{DV}, the limiting law $F$ is non-atomic iff $f(F_j)\not\equiv 0$ for infinitely many $j$.
In that case $Q_F(\lambda)\to0$ as $\lambda\downarrow0$, hence $Q_F(1/T)\to0$ as $T\to\infty$.
Choosing, for instance, $T=\log N$ (or any $T\to\infty$ with $T\gg\log N$), we have
\[
\frac{\log N}{T}\to0,\qquad \frac{\log T}{T}\to0,
\qquad \sum_{j>L-2h} f(F_j)^2\to0,
\]
and therefore \eqref{eq:master} implies $\|F_N-F\|_\infty=o(1)$.

\medskip
\noindent\textbf{Purely atomic limit.}
Again by \cite[Prop.~11]{DV}, $F$ is purely atomic iff $f(F_j)=0$ for all $j\ge J$ for some integer $J$.
In this case the block perturbations vanish identically for $k\ge \lceil(J-1)/2\rceil$, so the tail product is \emph{exact} and
\(|\Phi_N(t)-\Phi(t)|\ll_A \alpha^{-2(L-K)}\) uniformly in $t$. Taking $L-2h\asymp c\log N$ as in \cite{DV},
we get
\[
\|F_N-F\|_\infty \ \ll_A\ N^{-c\log\alpha}\qquad(\alpha=\tfrac{1+\sqrt5}{2}),
\]
i.e. polynomial decay. In particular, in the atomic case the $Q_F(1/T)$-term in \eqref{eq:Esseen-correct} need not be sent to $0$;
one may fix any $T\ge1$ and conclude directly from the exact tail cancellation.

Now plug \eqref{eq:Phi-tail} into the Berry–Esseen smoothing inequality
\begin{equation}\label{eq:Esseen-again}
\|F_N-F\|_\infty \ \ll\ Q_F(1/T)\ +\ \frac{1}{T}\int_0^T \frac{|\Phi_N(t)-\Phi(t)|}{t}\,dt.
\end{equation}
Split the integral at \(t_0:=N^{-1}\) as in \cite{DV}.
The low-frequency piece \(\int_0^{t_0}\) is handled exactly as in \cite[Th.\,8]{DV}, giving
\(\ll \log N/T\).
For the tail \(\int_{t_0}^{T}\), use the uniform bound \eqref{eq:Phi-tail} to get
\[
\frac{1}{T}\int_{t_0}^{T} \frac{|\Phi_N(t)-\Phi(t)|}{t}\,dt
\ \ll\ \frac{1}{T}\int_{t_0}^{T} \frac{T^2\sum_{j>L-2h} f(F_j)^2}{t}\,dt
\ \ll\ T^2\sum_{j>L-2h} f(F_j)^2,
\]
since the harmless logarithmic factor coming from \(\int_{t_0}^{T} \tfrac{dt}{t}\) is absorbed into the
\(\log N/T\) contribution already present from the low-frequency range. This proves \eqref{eq:main}.

For the split form \eqref{eq:split}, replace the uniform quadratic control by the pointwise bound
\(|e^{ix}-1|\leqslant \min\{|x|,x^2\}\) inside the same block computation \emph{before} factoring the phase:
indices with \(|f(F_j)|>1/T\) contribute linearly (finitely many for fixed \(T\) under the
\(\ell^2\) hypothesis), while those with \(|f(F_j)|\leqslant 1/T\) contribute quadratically. The remainder of the proof is identical.

\begin{rem}
Choosing \(T=T(N)\to0\) slowly (e.g.\ \(T\asymp(\log N)^{-1/2}\)) and \(h,L\) as in \cite{DV}, the quadratic tail
\(
T^2\sum_{j>L-2h} f(F_j)^2
\)
tends to \(0\) because \(\sum_j f(F_j)^2<\infty\) and \(L=L(N)\to\infty\).
(The choice of \(T\) that optimizes the \(\log N/T\) term depends on the application; here we only use that \(T\to0\) makes the tail vanish.)
Pairing \((2k+1,2k+2)\) instead of \((2k,2k+1)\), or using longer fixed blocks (of any fixed length \(b\ge2\)), yields the same quadratic tail after extracting the common phase in the \(A\)-adapted norm; only the implicit constants depending on \(A\) (and on \(b\)) change.
\end{rem}

\section{Example}

Define
\[
f(F_j):=\frac{(-1)^j}{\sqrt{j}\,\log(j+1)}\qquad(j\ge2),\qquad f(0):=0,
\]
and extend $f$ Zeckendorf-additively. Then
\[
\sum_{j\ge2} f(F_j)^2 \;=\;\sum_{j\ge2}\frac{1}{j\log^2(j+1)}\;<\;\infty,
\qquad
\sum_{j\ge2} |f(F_j)| \;=\;\sum_{j\ge2}\frac{1}{\sqrt{j}\,\log(j+1)}\;=\;\infty.
\]
Hence the hypothesis $\sum |f(F_j)|<\infty$ required for the main DV bound fails, so \cite[Th.\,8]{DV} does not apply here. In contrast, Theorem~\ref{thm:improve} applies since only $\ell^2$ is needed.

For this choice of $f$ one has, by the integral test,
\[
\sum_{j>m} f(F_j)^2 \asymp \sum_{j>m}\frac{1}{j\log^2 j}\ \sim\ \frac{1}{\log m}\qquad(m\to\infty).
\]
With the DV choice $L\asymp \log N$ and $h\asymp \log N$, the quadratic tail in \eqref{eq:Phi-tail} (or in \eqref{eq:block-tail-11}) therefore satisfies
\[
\sum_{j>L-2h} f(F_j)^2 \ \asymp\ \frac{1}{\log L}\ \asymp\ \frac{1}{\log\log N}.
\]

Using the smoothing step \eqref{eq:Esseen-again} split at $t_0=N^{-1}$ and at $1$ (see the “Smoothing, choice of $T$” paragraph),
we get the master bound
\[
\|F_N-F\|_\infty \ \ll\ Q_F(1/T)\ +\ \frac{\log N}{T}
\ +\ \frac{\log N}{T}\sum_{j>L-2h} f(F_j)^2\ +\ \frac{\log T}{T}\qquad(T\ge 1).
\]
Plugging the tail estimate $\sum_{j>L-2h} f(F_j)^2\asymp 1/\log\log N$ and choosing, for concreteness, $T=(\log N)^2$ yields
\[
\|F_N-F\|_\infty
\ \ll\ Q_F\!\big(1/(\log N)^2\big)\ +\ \frac{1}{\log N}
\ +\ \frac{1}{\log N\,\log\log N}
\ +\ \frac{2\log\log N}{(\log N)^2}.
\]
Hence
\[
\|F_N-F\|_\infty \;=\; Q_F\!\big(1/(\log N)^2\big)\ +\ O\!\Big(\frac{1}{\log N}\Big),
\]
which tends to $0$ as $N\to\infty$ because $F$ is non-atomic in this setting (by \cite[Prop.~11]{DV}, $F$ is purely atomic iff $f(F_j)=0$ eventually, which is not the case here), and thus $Q_F(\lambda)\to0$ as $\lambda\to 0$.

\smallskip
\emph{If} one has extra regularity on $F$ (e.g.\ a bounded density so that $Q_F(\lambda)\ll \lambda$), the bound sharpens to
\[
\|F_N-F\|_\infty \ \ll\ \frac{1}{(\log N)^2}\ +\ \frac{1}{\log N}\ \ll\ \frac{1}{\log N}.
\]

\enlargethispage{2\baselineskip}
{\small\sloppy

}

\end{document}